\documentclass[12pt]{article}
\usepackage{amscd, amssymb, latexsym, amsmath}

\newtheorem{thm}{Theorem}

\newtheorem{lem}{Lemma}

\newenvironment{proof}
{\noindent {\em Proof.}} {\hfill $\Box$}

\numberwithin{thm}{section} \numberwithin{cor}{section}
\numberwithin{pro}{section} \numberwithin{lem}{section}
\numberwithin{dfn}{section}
\numberwithin{rem}{section} \numberwithin{equation}{section}

\newcommand{\R}{\mathbb R}

\newcommand{\heat}{(\frac{d}{dt}-\Delta)}
\begin{document}
\title{A convergence result of the Lagrangian mean curvature flow}
\author{Mu-Tao Wang \footnote{The author is partially supported by National Science
Foundation Grant DMS0104163 and DMS0306049 and an Alfred P. Sloan
Research Fellowship.}}

\date{Aug. 15, 2005}
\maketitle \centerline{email:\, mtwang@math.columbia.edu }
\begin{abstract}

We prove the mean curvature flow of  the graph of a
symplectomorphism between Riemann surfaces converges smoothly as
time approaches infinity.

\end{abstract}

\section{Introduction}

Let $\Sigma_1$ and $\Sigma_2$ be two homeomorphic compact Riemann
surfaces without boundary. We assume $\Sigma_1$ and $\Sigma_2$ are
both equipped with Rimannian metrics of the same constant
curvature $c$, $c=-1, 0,$ or $1$. Let $\omega_1$ and $\omega_2$
denote the volume or symplectic forms of $\Sigma_1$ and
$\Sigma_2$, respectively. The Riemannian product space $\Sigma_1
\times \Sigma_2$ is denoted by $M$. We take
$\omega'=\omega_1-\omega_2$ to be the K\"ahler form of $M$ and $M$
becomes a K\"ahler-Einstein manifold with the Ricci form
$Ric=c\omega'$. Let $\Sigma$ be the graph of a symplectomorphism
$f:\Sigma_1\rightarrow \Sigma_2$, i.e. $f^*\omega_2=\omega_1$.
$\Sigma$ can be considered as a Lagrangian submanifold with
respect to the symplectic form $\omega'$.

The mean curvature flow deforms the initial surface
$\Sigma^0=\Sigma$ in the direction of its mean curvature vector.
Denote by $\Sigma^t$ the time slice of the flow at $t$. That
$\Sigma^t$ remains a Lagrangian submanifold follows from a result
of Smoczyk \cite{sm1}. The long-time existence and convergence
problems of this flow were studied in \cite{sm2} and \cite{mu2}.

In \cite{mu2}, the author proved the long time existence of the
flow and showed that $\Sigma^t$ for $t>0$ remains the graph of a
symplectomorphism $f_t$. When $c=1$, the author proved the
$C^\infty$ convergence as $t\rightarrow \infty$. However, only
$C^0$ convergence was achieved in the case when $c=-1$ or $0$.

 Independently, in \cite{sm2}, Smoczyk studied the case when $c=-1$
or $0$ assuming an extra angle condition. He discovered a
curvature estimate and showed that the second fundamental form is
uniformly bounded under this condition, and thus established the
long time existence and $C^\infty$ convergence at infinity.

In view of the above results, it is interesting to see whether the
$C^\infty $ convergence of the flow does require the angle
condition. In this paper, we show this assumption is unnecessary.

\begin{thm}
Let $(\Sigma_1, \omega_1)$ and $(\Sigma_2,\omega_2)$ be two
homeomorphic compact Riemann surface of the same constant
curvature $c=-1, 0,$ or $1$. Suppose $\Sigma$ is the graph of a
symplectomorphism $f:\Sigma_1\rightarrow \Sigma_2$ as a Lagrangian
submanifold of $M=(\Sigma_1\times \Sigma_2, \omega_1-\omega_2)$
and $\Sigma^t$ is the mean curvature flow with initial surface
$\Sigma^0=\Sigma$. Then $\Sigma^t$ remains the graph of a
symplectomorphism $f_t$ along the mean curvature flow. The flow
exists smoothly for all time and $\Sigma^t$ converges smoothly to
a minimal Lagrangian submanifold as $t\rightarrow \infty$.
\end{thm}

The long time existence part was already proved in \cite{mu2}. The
smooth convergence was established through a new integral estimate
(Lemma 3.1) related to the second variation formula. This estimate
is most useful when $c=-1 \,\text{or} \,0$. We remark the
existence of such minimal Lagrangian submanifold was proved using
variational method by Schoen \cite{sc} (see also Lee \cite{le}).

The author would like to thank Tom Ilmanen and Andre Neves for
useful discussions.
\section{Background material}
First we recall some formulas from \cite{mu2}.  The restriction of
the K\"ahler form $\omega'$ to $\Sigma^t$ gives a time-dependent
function $\eta=*\omega'$. Since $\Sigma^t$ is Lagrangian,
$*\omega'=2*\omega_1 $. $*\omega_1$ is indeed the Jacobian of the
projection $\pi_1$ from $M$ to $\Sigma_1$ when restricted to
$\Sigma$ and $\eta>0$ if and only if $\Sigma$ is locally a graph
over $\Sigma_1$. $\eta$ satisfies the following evolution
equation:

\begin{equation}\label{eta}
\frac{d}{dt}\eta=\Delta \eta +\eta[2|A|^2-|H|^2]+c\eta(1-\eta^2)
\end{equation} along the mean curvature flow.

Notice that $0< \eta \leq 1$. By the equation of $\eta$ and the
comparison theorem for parabolic equations, we get

\begin{equation}\label{lowerbound}\eta(x,t)\geq \frac{\alpha e^{ct}} {\sqrt{1+\alpha^2
e^{2ct}}}\end{equation} where $\alpha>0$ is given by
$\frac{\alpha} {\sqrt{1+\alpha^2}} =\min_{\Sigma_0}\eta$.
Therefore $\eta(x,t)$ converges uniformly to $1$ when $c=1$ and is
nondecreasing when $c=0$. In any case, $\eta$ has a positive lower
bound at any finite time and thus $\Sigma_t$ remains the graph of
a symplectomorphism.

Using the fact that the second fundamental form for Lagrangian
submanifold is a fully symmetric three tensor, one derives

\[|H|^2\leq \frac{4}{3}|A|^2.\]  Plug this into (\ref{eta}) and we
obtain

\begin{equation}
\frac{d}{dt}\eta\geq \Delta \eta
+\frac{2}{3}|A|^2\eta+c\eta(1-\eta^2)
\end{equation}

In \cite{mu2}, we apply blow-up analysis to this equation to show
there exists a weak blow-up limit with vanishing $\int |A|^2$.
This together with the lower bound of $\eta$ shows the limit is a
flat space and White's regularity theorem \cite{wh} implies the
blow-up center is a regular point. This proves the long-time
existence of the flow.

\section{A monotonicity lemma}

In this section, we derive a new monotonicity formula. First
$|H|^2$ satisfies the following evolution equation:

\begin{equation}\label{H^2}\heat |H|^2=-2|\nabla H|^2+2\sum_{ij}(\sum_k H_k h_{kij})^2
+c(2-\eta^2)|H|^2\end{equation} where the symmetric three-tensor
$h_{ijk}$ is the second fundamental form and $H_k=h_{iik}$, the
trace of the second fundamental form, is the component of the mean
curvature vector after identifying the tangent bundle and the
normal bundle through $J$.

We remark that both equations (\ref{eta}) and (\ref{H^2}) are
derived in Lemma 5.3 of \cite{sm2} where $p$ in \cite{sm2} and
$\eta$ are related by $\eta^2=\frac{4}{p}$ and $S=2c$.

We claim the following differential inequality is true:

\begin{lem}
\[\frac{d}{dt}\int_{\Sigma_t} \frac{|H|^2}{\eta}\leq c \int_{\Sigma_t} \frac{|H|^2}{\eta}\]
\end{lem}
\begin{proof}

The proof is a direct computation by combining equations
(\ref{eta}) and (\ref{H^2}). We compute

\[\begin{split}\frac{d}{dt}\frac{|H|^2}{\eta}&=\frac{\eta \Delta
|H|^2-|H|^2\Delta \eta}{\eta^2} -2\frac{|\nabla
H|^2}{\eta}\\
&+\frac{2\sum_{ij}(\sum_k H_k
h_{kij})^2-2|H|^2|A|^2+|H|^4}{\eta}+c\frac{|H|^2}{\eta}.\end{split}\]

Now

\[\Delta \frac{|H|^2}{\eta}=\frac{\eta\Delta|H|^2-
|H|^2\Delta \eta}{\eta^2}-\frac{2\eta\nabla\eta(\eta
\nabla|H|^2-|H|^2\nabla \eta )}{\eta^4}.\]

We plug this into the previous equation and obtain

\[\begin{split}\frac{d}{dt}\frac{|H|^2}{\eta}&=\Delta \frac{|H|^2}{\eta}
+\frac{2\eta\nabla\eta(\eta \nabla|H|^2-|H|^2\nabla \eta
)}{\eta^4} -2\frac{|\nabla
H|^2}{\eta}\\
&+\frac{2\sum_{ij}(\sum_k H_k
h_{kij})^2-2|H|^2|A|^2+|H|^4}{\eta}+c\frac{|H|^2}{\eta}.\end{split}\]

Rearranging terms, we arrive at

\[\begin{split}\frac{d}{dt}\frac{|H|^2}{\eta}&=\Delta \frac{|H|^2}{\eta}
+\frac{4\eta|H|\nabla\eta \cdot \nabla|H|-2|\nabla
\eta|^2|H|^2-2\eta^2|\nabla H|^2}{\eta^3} \\
&+\frac{2\sum_{ij}(\sum_k H_k
h_{kij})^2-2|H|^2|A|^2+|H|^4}{\eta}+c\frac{|H|^2}{\eta}.\end{split}\]

Integrate this identity and we have
\[\begin{split}\frac{d}{dt}\int_{\Sigma_t} \frac{|H|^2}{\eta}&=\int_{\Sigma_t}\frac{4\eta |H|\nabla\eta \cdot \nabla|H|-2|\nabla
\eta|^2|H|^2-2\eta^2|\nabla H|^2}{\eta^3} \\
&+\int_{\Sigma_t}\frac{2\sum_{ij}(\sum_k H_k
h_{kij})^2-2|H|^2|A|^2}{\eta}+c\int_{\Sigma_t}\frac{|H|^2}{\eta}.\end{split}\]

We use $|\nabla |H||\leq |\nabla H|$ in the first summand on the
right hand side and complete the square
\[\frac{4\eta|H|\nabla\eta \cdot \nabla|H|-2|\nabla
\eta|^2|H|^2-2\eta^2|\nabla |H||^2}{\eta^3}=-2\frac{|\nabla \eta
|H|-\eta\nabla |H||^2 }{\eta^3}.\] At last, we apply
Cauchy-Schwarz inequality to the second summand and the
differential inequality is proved.

\end{proof}
\section{Proof of the theorem}
The smooth convergence in the case when $c=1$, i.e. when
$\Sigma_1$ and $\Sigma_2$ are both standard $S^2$, was proved in
\cite{mu2}.

We prove the $C^\infty$ convergence in the case $c=0$ and $c=-1$
in the following. By the general convergence theorem of Simon
\cite{si}, it suffices to show $|A|^2$ is bounded independent of
time.

 In the case when $c=0$,
by (\ref{lowerbound}), $\eta$ has a positive lower bound. We have

\begin{equation}\label{H}\int_{\Sigma_t} |H|^2\leq
\int_{\Sigma_t}\frac{|H|^2}{\eta}\leq
K_1\int_{\Sigma_t}|H|^2\end{equation} for some constant $K_1$.

Since $\int_0^\infty \int_{\Sigma_t} |H|^2 <\infty$, there exists
a subsequence $t_i$ such that $\int_{\Sigma_{t_i}} |H|^2
\rightarrow 0$ and thus $\int_{\Sigma_{t_i}} \frac{|H|^2}{\eta}
\rightarrow 0$ as well. Because $\int_{\Sigma_t}
\frac{|H|^2}{\eta}$ is non-increasing, this implies
$\int_{\Sigma_{t}} \frac{|H|^2}{\eta} \rightarrow 0$ for the
continuous parameter $t$ as it approaches $\infty$.  Together with
(\ref{H}) this implies $\int_{\Sigma_{t}}{|H|^2} \rightarrow 0$ as
$t\rightarrow \infty$. By the Gauss formula,
$\int_{\Sigma_{t}}{|A|^2}=\int_{\Sigma_t}|H|^2 \rightarrow 0$. The
$\epsilon$ regularity theorem in \cite{il} (see also \cite{ec})
implies $\sup_{\Sigma_t} |A|^2$ is uniformly bounded.

In the case when $c=-1$, we have

\[\frac{d}{dt}\int_{\Sigma_t} \frac{|H|^2}{\eta}\leq -\int_{\Sigma_t}\frac{|H|^2}{\eta}
\]
or

\[\int_{\Sigma_t} \frac{|H|^2}{\eta}\leq K_2 e^{-t}\] for some
constant $K_2$,

 Since $\eta\leq 1$, we have
\[\int_{\Sigma_t} {|H|^2}\leq K_2 e^{-t}.\]

This implies $\Sigma_t \rightarrow \Sigma_\infty$ in Radon
measure. Indeed, for any function $\phi$ on $M$ with compact
support, it is easy to see

\[\frac{d}{dt}\int_{\Sigma_t} \phi=\int_{\Sigma_t}\nabla^M \phi
\cdot H+\int_{\Sigma_t} \phi |H|^2, \] and thus

\[\int_{\Sigma_t} \phi\rightarrow \int_{\Sigma_\infty} \phi.\]
exponentially. Also  the limit measure $\Sigma_\infty$ is unique.

The argument in \cite{mu2} shows the limit $\Sigma_\infty$ is
smooth. It seems one can adapt the proof of the local regularity
theorem of Ecker (Theorem 5.3) \cite{ec} or the original local
regularity theorem of Brakke \cite{br} to get the uniform bound on
second fundamental form. This does require  versions of these
theorem in a general ambient Riemannian manifold.

We circumvent this step by quoting a theorem in minimal surfaces.
Suppose the second fundamental form is unbounded. The blow-up
procedure in Proposition 3.1 of \cite{mu2} produces a limiting
flow that exists on $(-\infty, \infty)$. The flow has uniformly
bounded second fundamental form $A(x,t)$ and $|A|(0,0)=1$. It is
not hard to see each slice is the graph of an area-preserving map
from $\R^2$ to $\R^2$. Since $\int_{\Sigma_t} |H|^2\leq K_2
e^{-t}$, the limiting flow will satisfies

\[\int |H|^2\equiv 0\]
Therefore, we obtain a minimal area-preserving map. A result of Ni
\cite{ni} generalizing Schoen's theorem \cite{sc} shows this is a
linear diffeomorphism. This contradicts to the fact that
$|A|(0,0)=1$.


\begin{thebibliography}{99}

\bibitem{br} K. A. Brakke, \textit{ The motion of a surface by its mean
curvature.} Mathematical Notes, 20. Princeton University Press,
Princeton, N.J., 1978.

\bibitem{ec} K. Ecker, \textit{Regularity theory for mean curvature
flow.}
Progress in Nonlinear Differential Equations and their
Applications, 57. Birkhäuser Boston, Inc., Boston, MA, 2004.

\bibitem{es2}C. J. Earle and J. Eells, \textit{The diffeomorphism group of a compact
Riemann surface.} Bull. Amer. Math. Soc. 73 1967 557--559.














\bibitem{il} T. Ilmanen, \textit{Singularities
of mean curvature flow of surfaces.} preprint, 1997.



\bibitem{le} Y.-I. Lee, \textit{Lagrangian minimal surfaces
in Kähler-Einstein surfaces of negative scalar curvature.} Comm.
Anal. Geom. 2 (1994), no. 4, 579--592.









\bibitem{ni} L.  Ni, \textit{A Bernstein type theorem for minimal volume preserving
maps.} Proc. Amer. Math. Soc. 130 (2002), no. 4, 1207--1210.

\bibitem{sc} R. Schoen, \textit{
The role of harmonic mappings in rigidity and
 deformation problems.} Complex geometry
 (Osaka, 1990), 179--200, Lecture Notes in
 Pure and Appl. Math., 143, Dekker, New York,
 1993.

\bibitem{si} L. Simon, \textit{Asymptotics for
 a class of nonlinear evolution equations,
 with applications to geometric problems.}
  Ann. of Math. (2) 118 (1983), no. 3, 525--571.



\bibitem{sm1} K. Smoczyk, \textit{Der Lagrangesche mittlere Kr\"ummungsflu\ss\ (The
Lagrangian mean curvature flow).} Habilitation thesis (English
with German preface), University of Leipzig, Germany (1999), 102
pages.


\bibitem{sm2} K. Smoczyk, \textit{Angle theorems for the Lagrangian mean curvature
flow.} Math. Z. 240 (2002), no. 4, 849--883.




\bibitem{mu1}  M.-T. Wang, \textit{Mean curvature
flow of surfaces in Einstein four-Manifolds.} J. Differential
Geom. \textbf{57} (2001), no. 2, 301-338.

\bibitem{mu2} M.-T. Wang, \textit{Deforming area preserving
diffeomorphism of surfaces by mean curvature flow.} Math. Res.
Lett. 8 (2001), no.5-6, 651-662.






\bibitem{wh} B. White, \textit{A local
regularity theorem for classical mean curvature flow.} preprint,
1999, revised 2002.
\end{thebibliography}
\end{document}